\documentclass[12pt]{amsart}

\newcommand{\RNum}[1]{\uppercase\expandafter{\romannumeral #1\relax}}
\usepackage{yfonts}
\newtheorem*{conjecture}{Conjecture}
\newtheorem{lemma}{Lemma}
\newtheorem*{note}{Note}

\usepackage{pdfpages}
\usepackage{amssymb,latexsym}
\usepackage{listings}
\usepackage{fancyvrb}
\newtheorem{theorem}{ Theorem}
\newtheorem{definition}{Definition}

\begin{document}

\begin{center}
{\Large A study of a family of generating functions of Nelsen-Schmidt type and some identities on restricted barred preferential arrangements}\\ 
 \vspace{5mm}
{\large S.Nkonkobe, \, V.Murali}\\
\vspace{2mm}\

{\it \footnotesize Department of Mathematics (Pure \& Applied)\\ Rhodes
University \\Grahamstown 6140 South Africa\\ snkonkobe@yahoo.com, v.murali@ru.ac.za}    \\

\end{center}
\section{Abstract\label{section:1}}
A preferential arrangement of a set $X_n=\{1,2,\ldots,n\}$ is an ordered partition of the set $X_n$ induced with a linear order.  Separation of blocks of a preferential arrangement with bars result in the notation of barred preferential arrangements. Roger Nelsen and Harvey Schmidt proposed the family of generating functions $P^k(m)=\frac{e^{km}}{2-e^m}$; which for $k=0$ and for $k=2$ they have shown that the generating functions are exponential generating functions for the number of preferential arrangements of a set $X_n$ and the number of chains in the power set of $X_n$ respectively. In this study we propose  combinatorial structures whose integer sequences are generated by members of the family for all values of $k$ in $\mathbb{Z}^+$. To do this we use a notion of restricted barred preferential arrangements. We then propose a more general family of generating functions 
 $P^{r}_{j}(m)=\frac{e^{rm}}{(2-e^m)^j}$ for $r,j\in\mathbb{Z^+}$. We derive some new identities on restricted barred preferential arrangements and give their combinatorial proofs. We also propose conjectures on number of restricted barred preferential arrangements.
  
\fontsize{8}{1}\selectfont     
Mathematics Subject Classifications:05A18,05A19,05A16, 2013
\\\textbf{Keywords}: Barred preferential arrangements, restricted barred preferential arrangements, restricted sections and free sections.
\normalsize
\section{introduction\label{section:2}}
The integer sequence 1,1,3,13,75,541,... that arise for the number of preferential arrangements of $X_n$ can be tracked back to \cite{cayley:1859} a year 1859 paper by Arthur Cayley; in connection with analytical forms called trees. The word preferential arrangements is due to Gross in \cite{gross:1962}. The integer sequence is A000670 in Sloane \cite{Sloane:preferential}. Since Cayley's paper much work has been done: The sequence has been given various interpretations and has been connected to a number of well known combinatorial sequences (for instance in \cite{mendelson:1982}, \cite{murali:Nelsongenerating function}, \cite{Pippenger:2010}, \cite{barred:2013}). Examples of preferential arrangements of $X_3=\{1,2,3\}$ are:\\ a)\; 3\quad 12\\ b)\:\;2\quad3\quad1

 The preferential arrangement in a) constitutes of two blocks. The first block of the preferential arrangement (from left to right) is the element 3. The second block of the preferential arrangement is formed by the elements 1 and 2. The preferential arrangement in b) has three blocks. The first block being the element 2. The second block being the element 3 and the third block is the element 1. The magnitude of the spaces between the elements tell us which elements form a single block.

In \cite{Nelsen:91} Nelsen and Schmidt have proposed the family of generating functions $P^k(m)=\frac{e^{km}}{2-e^m}$; which for $k=0$ and for $k=2$ they have shown that the generating functions are exponential generating functions for the number of preferential arrangements on $X_n$ and the number of chains in the power set of $X_n$ respectively. They then asked ``could there be other combinatorial structures associated with either the set $X_n$ or the power set of $X_n$ whose integer sequences are generated by  members of the family for other values of $k$ ?"
\\\\Murali has shown that for four values of $k$ (i.e. $k=0,1,2,3$) the members of the family are generating functions for the number of restricted preferential fuzzy subsets of a set $X_n$ \cite{murali:Nelsongenerating function}. The main aim of this study is to  propose combinatorial structures whose integer sequences are  generated by the members of the family for all values of $k$ in $\mathbb{Z^+}$. From now on for a fixed $k\in\mathbb{Z^+}$ we will refer to the generating function $P^k(m)=\frac{e^{km}}{2-e^m}$ of Nelsen and Schmidt, as the Nelsen-Schmidt generating function.  
\section{Preliminaries\label{section:3}}
Separating blocks of a preferential arrangement with a number of bars results in a barred preferential arrangement \cite{barred:2013}.  Examples of barred preferential arrangements of the set $X_3$ having one bar and two bars are:$\\a)\: 2|\:1\:\:3\\b)\:3|\:|\:12$

In a) the three elements of $X_3$ are all forming separate blocks. There is a single bar separating the first and the second block. In b) the set $X_3$ is preferentially arranged into two blocks. One block is the element 3 and the other block is formed by the elements 1 and 2. The two bars are in-between the two blocks. With reference to bars the barred preferential arrangement in a)  has  two sections one to the left of the bar (the block 2) and another section to the right of the bar (the two blocks). We name them section zero and section one respectively.
 The barred preferential arrangement in b) has three sections. One section to the left of the first bar. A section in-between the two bars and another section to the right of the second bar. We name them section zero, section one and section two respectively. In general $k$ bars separate a preferential arrangement into $k+1$ sections (possibly with empty sections) \cite{barred:2013}. 
\begin{note}\nonumber The elements within each section of a preferential arrangement are preferentially arranged among themselves.
\end{note}

 We denote the set of all barred preferential arrangements of the set $X_n$ having $k$ bars by $Q^k_n$ and the number of these barred preferential arrangements by $J^k_n$ i.e $|Q^k_n|=J^k_n$. A closed form, a recurrence relation and generating function for the numbers $J^k_n$ are given respectivel in the following theorems; \begin{theorem}\cite{barred:2013}\label{theorem:1} for all $n,k\geq0$ we have; \begin{center}$J^k_n=\sum\limits_{s=0}^n {n \brace s}s!\binom{k+s}{s}$\end{center}\end{theorem}
 Where ${n \brace s}$ are Stirling numbers of the second kind.
 \begin{theorem}\cite{barred:2013}\label{theorem:2} for $n,k\geq0$ we have;\begin{center} $J^k_n=\sum\limits_{s=0}^{n}\binom{n}{s}J^0_{s}J^{k-1}_{n-s}$\end{center}
 
 \end{theorem}\begin{theorem}\cite{barred:2013}\label{theorem:3} for $k\geq0$ we have;  \begin{center}$q^k(m)=\sum\limits_{n=0}^{\infty}\frac{J_n^k m^n}{n!}=\frac{1}{(2-e^m)^{k+1}}$\end{center}\end{theorem}
 \section{Restricted barred preferential arrangements\label{section:4}}
 Whether we view  barred preferential arrangements as a result of first placing bars and then distributing elements on each section or first preferentially arranging elements into blocks and then introduce bars in-between the blocks of the preferential arrangements; the resulting number of barred preferential arrangements is the same. Here we choose the former.
 
 In finding the total number of barred preferential arrangements of $X_n$ having $k$ bars we generalise equation 24 of \cite{Pippenger:2010} from one bar to $k$ bars using the same kind of argument as in \cite{Pippenger:2010}; we argue as follows:
  We first place the $k$ bars. There are $k+1$ sections before, in-between and after the $k$ bars. In a distribution of the $n$ elements of $X_n$, say there are $w_i$ elements on the $i^{th}$ section. There are $n!$ ways of permuting the $n$ elements among the $k+1$ sections. There are $w_i!$ ways of permuting elements of section $i$ among themselves. There are $J^0_{w_i}$ ways of preferentially arranging elements of section~$i$ among themselves. Hence the total number $J^k_n$ of barred preferential arrangements of $X_n$ having $k$ bars is;\begin{equation}\label{equation:1} J^k_n=\sum\limits_{w_{1}+\cdots+w_{k+1}= n}\frac{n!}{w_{1}!\times w_{2}!\cdots w_{k+1}!}J^0_{w_{1}}\times J^0_{w_{2}}\times\cdots\times J^0_{w_{k+1}}\end{equation}   
Where the summation is taken on all solutions of the equation \\$w_{1}+\cdots+w_{k+1}= n$ in non-negative integers.

Equation \ref{equation:1} is a dual of the closed form in theorem \ref{theorem:1} in section \ref{section:3} above. In (\ref{equation:1}) the number of barred preferential arrangements is obtained using an argument where bars are first placed and then elements being distributed on the  sections. Whereas in theorem \ref{theorem:1} the number of barred preferential arrangements is obtained using an argument; where the elements of $X_n$ are first  preferentially arranged into blocks, then bars are introduced in-between the blocks.

\noindent { {$ \mathbf 1^ \circ $ } }
\begin{definition}\label{definition:1}
A free section of a barred preferential arrangement is a section whose elements can be preferentially arranged among themselves in any possible way. 
\end{definition} 
What definition \ref{definition:1} means is that; when elements of a free section are preferentially arranged, the preferential arrangement can have any possible number of blocks. So for $w_i$ elements in a free section there are $J^0_{w_i}$ possible ways of preferentially arranging the elements among themselves. 
\begin{definition}\label{definition:2}
A restricted section of a barred preferential arrangement is a section that can have a maximum of one block.
\end{definition}
What definition \ref{definition:2} means is that; when elements of a restricted section are preferentially arranged; they can have a maximum of one block. 
So for $w_j$ elements on a restricted section; there is 1 way of preferentially arranging the elements among themselves i.e. into one  block.
\\\\In constructing a barred preferential arrangement of $X_n$ having $k$ bars by first placing bars and then distributing elements on sections; we require that one fixed section to be a free section and all the other $k$ sections to be restricted sections. The chosen fixed section can be any of the sections but must be fixed. In getting the number of barred preferential arrangements with the restriction we argue as follows: Lets say there are $w_i$ elements on section $i$. We assume the free section is the $j^{th}$ section (from left to right). There are $n!$ ways of permuting the $n$ elements of $X_n$ among the $k+1$ sections. There are $w_i!$ ways of permuting elements of section $i$ among themselves. There are $J^0_{w_j}$ ways of preferentially arranging elements of the chosen free section among themselves. There is one way of preferentially arranging elements of each of the $k$ restricted sections. Hence the total number of these restricted barred preferential arrangements (denoted by $p^k_n$) is given by;

\fontsize{8.5}{10}\selectfont \begin{equation}\label{equation:2}
p^k_n=\sum\limits_{w_{1}+\cdots+w_{k+1}= n}\frac{n!}{w_{1}!\times w_{2}!\times\cdots\times w_{j}!\times\cdots\times w_{k+1}!}\begin{bmatrix}(1)\times(1)\times\cdots \times J^0_{w_{j}}\times\cdots\times (1)\end{bmatrix}
\end{equation}\normalsize
Where the summation is taken on all solutions of the equation \\$w_{1}+\cdots+w_{k+1}= n$ in non-negative integers.
In the product of the $k+1$ terms in the bracket in (\ref{definition:2}); all the entries are ones except the $j^{th}$ entry which is $J^0_{w_{j}}$. This is due to the fact that the $j^{th}$ section is the chosen free section. The other $k$ sections are all restricted sections.
 We denote the set of these restricted barred preferential arrangements by $G^k_n$ so $|G^k_n|=p^k_n$. 
The set $G^k_n\subseteq Q^k_n$ where $Q^k_n$ is the set of all barred preferential arrangements of $X_n$ having $k$ bars without any restrictions considered in section \ref{section:3} above. 

We now seek the exponential generating function for the numbers $p^k_n$ in (\ref{equation:2}). In finding the generating function of the numbers we define a convolution of exponential generating functions: Given generating functions $X_1=\sum\limits_{n_1=0}^{\infty}\frac{x_{n_1}m^{n_1}}{n_1!}$,\ldots, $X_s=\sum\limits_{n_s=0}^{\infty}\frac{x_{n_s} m^{n_s}}{n_s!}$ their convolution is;

\fontsize{9}{9}\selectfont
\begin{equation}\label{equation:3}  X_1\times X_2\times\cdots\times X_s=\sum\limits_{n=0}^{\infty}\begin{pmatrix}\sum\limits_{n_1+\cdots+n_s= n}\frac{n!}{n_1!\times n_2!\cdots n_s!}x_{n_1}^{1}\times x_{n_2}^{2}\cdots x_{n_s}^{s}\end{pmatrix}\frac{m^n}{n!}\end{equation}\normalsize
Where the summation is taken on all solutions of the equation \\$n_1+\cdots+n_s= n$ in non-negative integers.
By (\ref{equation:3}) the generating function of the sequence $p^k_n$ in (\ref{equation:2}) is;

\fontsize{5}{6}\selectfont
\begin{equation}\label{equation:5}
P^k(m)=\begin{pmatrix}\sum\limits_{n_1=0}^{\infty}\frac{ m^{n_1}}{n_1!}\end{pmatrix}\times\begin{pmatrix}\sum\limits_{n_2=0}^{\infty}\frac{ m^{n_2}}{n_2!}\end{pmatrix}\times\cdots\times\begin{pmatrix}\sum\limits_{n_j=0}^{\infty}\frac{ J^0_{w_j} m^{n_j}}{n_j!}\end{pmatrix}\times\cdots\times\begin{pmatrix}\sum\limits_{n_{k+1}=0}^{\infty}\frac{ m^{n_{k+1}}}{n_{k+1}!}\end{pmatrix}=\sum\limits_{n=0}^{\infty}\frac{p^k_n m^n}{n!}
\end{equation}\normalsize 
We observe that each term in the product in (\ref{equation:5}) is a generating function. Where $k$ of them are $e^m$ and one of them is $\frac{1}{2-e^m}$.
 Hence the generating function in (\ref{equation:5}) is;
 
 \fontsize{7}{6}\selectfont
\begin{equation}\label{equation:4}
P^{k}(m)=\sum\limits_{n=0}^{\infty}\frac{p^k_n\: m^n}{n!}=e^m\times e^m\times\cdots\times\frac{1}{2-e^m}\times\cdots\times e^m=\frac{e^{km}}{2-e^m}
\end{equation}\normalsize
The term $\frac{1}{2-e^m}$ in the generating function in (\ref{equation:4}) is due to the free section (see theorem \ref{theorem:3} of section \ref{section:2} above: this is the case $k=0$). The terms $e^m$ are due to the restricted sections. 

We recognise the generating function in (\ref{equation:4}) as the  Nelsen-Schmidt generating function. So the Nelsen-Schmidt generating function for a fixed $k$ in $\mathbb{Z^+}$; is the generating function for the number of barred preferential arrangements of $X_n$ having $k$ bars. In which one fixed section is a free section and the other $k$ sections are restricted sections.
The idea of restricted barred preferential arrangements with one bar, where one section is a free section and the other is a restricted section appear in \cite{benjamin:Nelsen problem for t equal 1} in a different context. The idea in the paper appears in a campus security problem, where the author seeks the total possible number $A_n$ of possible combination locks having $n$ button. In the paper they derived the exponential generating function of $A_n$ as $\frac{e^x}{2-e^x}$. Hence our work above can be viewed as a generalisation of the idea to multiple bars.

\begin{theorem}\label{theorem:4}for $n\geq0$, $k\geq1$ we have;
\begin{center}$p^k_n=\sum\limits_{s=0}^{n}\binom{n}{s}p^0_{s}k^{n-s}$\end{center}
\end{theorem}We recall $p^k_n$ is the number of barred preferential arrangements of $X_n$ having $k$  bars, in which one fixed section is a free section; and the other $k$ sections are restricted sections. The set of these restricted barred preferential arrangements being $G^k_n$.  
In proving the theorem we view the number of elements in $G^k_n$ being obtained in the following way: On each $\textgoth{B}\in G^k_n$ there are $\binom{n}{s}$ ways of selecting elements to go into the free section. There are $p^{0}_{s}$ ways of preferentially arranging the $s$ elements among themselves. Since the other $k$ sections other than the chosen free section can have a maximum of one block; then the number of possible ways of distributing the remaining $n-s$ elements among the $k$ restricted sections is $k^{n-s}$. Taking the product and summing over $s$ we obtain the result of the theorem.
\begin{theorem}\label{theorem:5} for $n,k\geq0$ we have; 
\begin{center}
$p^{k+1}_{n}=\sum\limits^{n}_{s=0}\binom{n}{s}p^k_{s}$
\end{center}
\end{theorem}
We recall $p^{k+1}_n$ is the number of barred preferential arrangements an $n$-element set having $k+1$ bars in-which $k+1$ sections are restricted sections and one section is a free section. The set of these barred preferential arrangements being denoted by $G^{k+1}_n$. The prove of the theorem is similar to that of theorem~\ref{theorem:4}. In proving the theorem we base our argument on a fixed restricted section, say its the $m^{th}$ section on each element of $G^{k+1}_n$. Lets say there are $s$ elements which are not to go to the$m^{th}$ section on elements of $G^{k+1}_n$. The $s$ elements can be chosen in $\binom{n}{s}$ number of ways. The $s$ elements can be preferentially arranged on the $k+1$ other sections in $p^k_n$ ways. The remaining $n-s$ elements can be preferentially arranged on the $m^{th}$ section in one way. Taking the product and summing over $s$ we obtain the result of the theorem.

\begin{lemma}\label{lemma:1}for $k,n\geq1$ we have
\begin{center}$p^k_n=k^n+\sum\limits_{s=0}^{n-1}\binom{n}{s}k^s p^o_{n-s}$\end{center}
\end{lemma}
We recall $p^k_n$ is the number of restricted barred preferential arrangements of an $n$ element set having $k$ bars in-which $k$ fixed sections are restricted sections and one fixed section is a free section. Where the set of these barred preferential arrangements is denoted by $G^k_n$.  In proving the lemma we base our argument on a fixed free section on all barred preferential arrangements in $G^k_n$. For argument sake lets say the chosen free section is the first section on each element of $G^k_n$. On all barred preferential arrangements from $G^K_n$  either the first section is empty or non-empty. When the first section is empty that means the $n$ elements are preferentially arranged among the other $k$ sections of which $k$ all of them are restricted sections so the number of elements of $G^k_n$ in this case is $k^n$.
\\\\ The other case is when the first section has at least one element. In getting the number of elements of $G^k_n$ in this case we argue as follows; There can be a maximum of $n-1$ elements not in the first section in the this case and a minimum of 0 elements. We assume there are $s$ elements not in the first section in this case. There are $\binom{n}{s}$ of selecting elements which are not in the first section. There are $k^s$ ways of preferentially arranging the $s$ elements among the other $k$ sections.  The remaining $n-s$ elements can be preferentially arranged in the first section in $p^0_{n-s}$ ways. Taking the product and summing over $s$ we have the number of elements of $G^k_n$ in this case as $\sum\limits_{s=0}^{n-1}\binom{n}{s}k^sp^0_{n-s}$. \\\\Combining the two cases we obtain the result of the lemma.

\begin{lemma}\label{lemma:2} for $n,k\geq1$ we have
\begin{center}$p^k_n=p^{k-1}_n+\sum\limits_{s=0}^{n-1}\binom{n}{s}p^{k-1}_{s}$\end{center}
\end{lemma}
Lemma~\ref{lemma:2} can be proved in a similar way to lemma~\ref{lemma:1}.  We will generalise lemma~\ref{theorem:2} in \noindent {{$\mathbf 2^ \circ$}} as theorem~\ref{theorem:8} give a proof of the theorem. 

\begin{conjecture}\label{conjecture:1}for $k,s,n\geq0$ we have\begin{center}
$p^k_n=\frac{1}{2}\sum\limits_{s=0}^{\infty}\frac{(k+s)^n}{2^s}$\end{center}
\end{conjecture}
For the case $k=0$ the identity in the conjecture is equation 21 of~\cite{Pippenger:2010}.
\\\\\noindent { {$ \mathbf 2^ \circ $ } }\\
We ask: \emph{What sort of a generating function would restricted barred preferential arrangements having more than one free section have? As opposed to the restricted barred preferential arrangements considered in\noindent { {$ \mathbf 1^ \circ $}} above. Also how would the generating function for a fixed number of restricted sections; and fixed number of free sections compare to the Nelsen-Schmidt generating function we considered above?}
\\\\In constructing a barred preferential arrangement of $X_n$ having $k$ bars by first placing bars and then distributing elements on sections; from the $k+1$ sections we require that $j$ fixed sections to be free sections and remaining $r = k+1-j$ sections to be restricted sections. We denote the number of these restricted barred preferential arrangements by $p^{r}_{j}(n)$ and the set of these restricted barred preferential arrangements by $G^r_j(n)$; so $|G^r_j(n)|=p^{r}_{j}(n)$. In a similar way to the way we obtained the number $p^k_n$ in (\ref{equation:2}); the number $p^{r}_{j}(n)$ of these restricted barred preferential arrangements is;

\fontsize{7.5}{1}\selectfont
\begin{equation}\label{equation:7}
p^{r}_{j}(n)=\sum\limits_{w_{1}+\cdots+w_{k+1}= n}\frac{n!}{w_{1}!\times w_{2}!\cdots w_{m}!\cdots w_{k+1}!}(1)\times(1)\times\cdots\times(1)\times J^0_{w_{m}}\times J^0_{w_{m+1}}\times\cdots\times J^0_{w_{k+1}}
\end{equation}\normalsize 
Where on the product of $k+1$ terms in (\ref{equation:7}); $r$ terms are (1)'s for the $r$ restricted sections and $j=k+1-r$ terms are the ${J^0_{w_{i}}}'s$ for the free sections.
 
In a similar way to the way the generating function for the numbers $p^k_n$ is obtained in (\ref{equation:4}); the generating function (denoted by $P^r_{j}(m)$) for the numbers $p^{r}_{j}(n)$ is;
\begin{equation}\label{equation:8}
P^r_{j}(m)=\sum\limits_{n=0}^{\infty}\frac{p_j^{r}(n)\:m^n}{n!}=\frac{e^{rm}}{(2-e^m)^{j}}
\end{equation}
We observe that Nelsen-Schmidt generating function is a special case of the generating function $P^r_{j}(m)$ when $j=1$.
\\\\We  generalise lemma \ref{lemma:1} of \noindent { {$ \mathbf 1^ \circ $}} above to the following theorem.  
\begin{theorem}\label{theorem:6} for $r,j,n\geq1$ we have
\begin{center}
$p^r_j(n)=p^{r}_{j-1}(n)+\sum\limits_{s=0}^{n-1}\binom{n}{s}p^{r}_{j-1}(s) p_{1}^{0}(n-s)$
\end{center}
\end{theorem}
We recall $p^r_j(n)$ is the number of restricted barred preferential arrangements of an $n$ element set having $r+j-1$ bars in-which $r$ fixed sections are restricted sections and $j$ fixed sections  are free sections. Where the set of these barred preferential arrangements is denoted by $G^r_j(n)$. In proving the theorem we base our argument on a fixed free section on all barred preferential arrangements in $G^r_j(n)$. For argument sake lets say the chosen free section is the first section of each element of $G^r_j(n)$ (the chosen fixed free section could have been any of the $j$ free sections as long would be fixed). On all barred preferential arrangements from $G^r_j(n)$  either the first section is empty or non-empty. So we have two cases.
\\\\When the first section is empty that means the $n$ elements are preferentially arranged among the other $r+j-1$ sections of which $r$ of them are restricted and $j-1$ of them are free sections. So the number of elements of $G^r_j(n)$ in this case is $p^r_{j-1}(n)$ (by definition of $p^r_{j-1}(n)$).
\\\\ The other case is when the first section of each element of $G^r_j(n)$  has at least one element. In getting the number of elements of $G^r_j$ in this case we argue as follows; There can be a maximum of $n-1$ elements not in the first section in the this case and a minimum of 0 elements. We assume there are $s$ elements not in the first section in this case. There are $\binom{n}{s}$ of selecting elements which are not in the first section. There are $p^{r}_{j-1}(n)$ ways of preferentially arranging the $s$ elements among the other $r=j-1$ sections.  The remaining $n-s$ elements can be preferentially arranged in the first section in $p^0_1(n-s)$ ways. Taking the product and summing over $s$ we have the number of elements of $G^r_j(n)$ in this case as $\sum\limits_{s=0}^{n-1}\binom{n}{s}p^{r}_{j-1}(n)p^0_1(n-s)$. \\\\Combining the two cases we obtain the result of the lemma.
\\\\We generalise lemma~\ref{lemma:2} of \noindent { {$ \mathbf 1^ \circ $}} to the following theorem. 

\begin{theorem}\label{theorem:8} for $n,r,j\geq1$ we have
\begin{center}
$p^r_j(n)=p^{r-1}_{j}(n)+\sum\limits_{s=0}^{n-1}\binom{n}{s}p^{r-1}_j(s)$
\end{center}
\end{theorem}
The number $p^r_j(n)$ is the number of barred preferential arrangements of an $n$-element set having $r+j-1$ bars in-which $r$ sections are restricted sections and $j$ sections are free sections. Where the set of these barred preferential arrangements is denoted by $G^r_j(n)$. The proof of this theorem is similar to that of theorem~\ref{theorem:6} above. We base our argument on a fixed restricted section. We say the chosen restricted section is the first section of each element of $G^r_j(n)$ (the chosen restricted section could have been any of the $r$ restricted sections as long it would be fixed). All barred preferential arrangements of the set $G^r_j(n)$ have the property that the chosen restricted section is either empty or has a certain  number of elements. We partition the set $G^r_j(n)$ into two disjoint subsets $w_1$  and $w_2$. Where the sets contain those elements of $G^r_j(n)$ in-which the chosen restricted section is empty and those in-which the restricted section has a number of elements respectively. In obtaining the cardinality of the set $w_1$ we argue as follows: On those elements of $G^r_j(n)$ in which the chosen restricted section is empty; that means the $n$ elements are distributed among the other $r+j-1$ sections of-which $r-1$ of the sections are restricted sections and $j$ of the sections are free sections.  Hence the cardinality of $w_1$ is $p^{r-1}_{j}(n)$ (by definition of $p^{r-1}_{j}(n)$). 

In finding the number of elements of the set $w_2$ we argue as follows; The maximum number of elements which are not in the chosen section in this case is $n-1$ and the minimum number is 0. Lets say there are $s$ elements which are not in the restricted section. There are $\binom{n}{s}$ ways of selecting the elements. There are $p^{r-1}_{j}(s)$ ways of preferentially arranging the $s$ elements among the $r+j-1$ sections. There are $p^1_0(n-s)$ ways of preferentially arranging the remaining $n-s$ elements on the restricted section (note $p^1_0(n-s)=1$ ). Taking the product and summing over $s$ we obtain  $\sum\limits_{s=0}^{n-1}\binom{n}{s}p^{r-1}_{j}(s)=w_2$.\\ Combining $w_1$ and $w_2$ we obtain the result of the theorem.

\begin{theorem}\label{theorem:7} for $n,r,\geq0$ and $j\geq1$ we have;\begin{center}$p^{r+1}_{j}(n) =2p^{r}_{j}(n)-p^r_{j-1}(n)$
\end{center}\end{theorem} 
In proving the theorem we consider the set $G^r_j(n)$; which is the set of all barred preferential arrangements of a set $X_n$ having $k$ bars in which $r$ fixed sections are restricted sections and $j$ fixed sections are free sections. By choice we choose the $j$ free sections (see definition \ref{definition:1} above) of each element $\textgoth{B}\in G^r_j(n)$ to be the last (from left to right) $j$ sections of $\textgoth{B}$. We want to construct the set $G^{r+1}_j(n)$ using elements of the set $G^r_j(n)$.

We argue as follows: We add an extra bar $\overset{*}{{|}}$ to the far right of each element $\textgoth{B}\in G^r_j(n)$ to form the set $D^r_j(n)$. We observe that adding the extra bar to each element of $G^r_j(n)$ does not affect counting; hence $|D ^r_j(n)|=|G^r_j(n)|=p^{r}_{j}(n)$.
 We note that on each element $\textgoth{C}\in D^r_j(n)$ to the left of the bar $\overset{*}{{|}}$ is a free section and the section to the right of $\overset{*}{{|}}$ is empty.
 We define the set $R^r_{j}(n)=\{0,1\}\times D^r_j(n)$ as containing the same elements as the set $D^r_j(n)$ but with an indexing on the bar $\overset{*}{{|}}$ which is 0 and separately 1. Hence $R^r_{j}(n)$ has twice the number of elements as $D^r_j(n)$ (of which half of them have the index 0 and the other half have the index 1). Also each element of $R^r_{j}(n)$ has $k+2$ sections due to the introduction of the bar $\overset{*}{{|}}$. 
We now use elements of $R^r_{j}(n)$ to construct the set $G^{r+1}_j(n)$. We construct as follows:
\\\RNum{1}.\; If the index on the bar $\overset{*}{{|}}$ on an element $\textgoth{O}\in R^r_{j}(n)$ is 0 then such an element will be interpreted in $G^{r+1}_j(n)$ as an element of  $G^{r+1}_j(n)$ whose ${(k+2)}^{th}$ section is empty.

 We collect all such elements to form the set $W$. The set $W$ has $p^{r}_{j}(n)$ elements (since half of the elements in $R^r_{j}(n)$ have the index 0)
 \\\RNum{2}. If the index on the bar $\overset{*}{{|}}$ on $\textgoth{O}\in R^r_{j}(n)$ is 1; then we shift the last block of the ${(k+1)}^{th}$ section of $\textgoth{O}$ to be the only block to the right of $\overset{*}{{|}}$ (i.e the block closest to the bar $\overset{*}{{|}}$). We collect all such elements to form the set $K$. There are $p^{r}_{j}(n)$ elements having index 1 in $R^r_{j}(n)$ (half of the elements in $R^r_{j}(n)$ have the indexing 1).\\\\ When discarding the index on the bar  $\overset{*}{{|}}$ in the construction of the set $K$; some elements are common elements between the set $K$ and the set $W$. The common elements occur when the free section closest to the bar  $\overset{*}{{|}}$ is empty. In that case there is no block to put to the right of the bar $\overset{*}{{|}}$ when constructing the set $K$. The common elements between $K$ and $W$ occur when the $n$ elements of $X_n$ are distributed on the first $k$ sections of each element of $G^r_j(n)$. The $k$ other sections are composed of $r$ restricted sections and $j-1$ free sections. So the number of ways in which the elements of $X_n$ are distributed on these sections is $p^{r}_{j-1}(n)$; so $|K\cap W|=p^{r}_{j-1}(n)$.
  Hence $|K\cup W|=|K|+|W|-|K\cap W|$$\Rightarrow$ $|K\cup W|=p^{r}_{j}(n)+p^{r}_{j}(n)-p^{r}_{j-1}(n)$. The elements of  both sets $K$ and $W$ have $k+2$ sections; Of which for a fixed $\textgoth{X}\in K\cup W$,  $r+1$ sections  are restricted sections and $j$ sections are free sections. By definition $K\cup W=G^{r+1}_{j}(n)$. This completes the proof. 
 \\\\ The statement of theorem \ref{theorem:7} for the case $j=1$ is $p^{r+1}_{1}(n)=2p^{r}_{1}-r^n$. Which is  equation 14 of \cite{Nelsen:91}, proposed by Nelsen and Schmidt.

\begin{conjecture}\nonumber for $n,r\geq0$ we have; \begin{center}$p_1^r(n)=\sum\limits^{n}_{s=1}\binom{n}{s}p_1^r(n-s)+r^n$\end{center}\end{conjecture} 
The identity in the conjecture when $r=0$ appears in \cite{gross:1962} as equation~9.  For $r=2$ the identity in the conjecture occurs in \cite{murali:combinatorics} as an identity for the number of preferential fuzzy subsets of a set $X_n$ ; it occurs as equation~3.4 in \cite{murali:combinatorics}. 

\section{Future work}

 MacMahon in \cite{MacMahon} has shown a way on how the generating function $P^k(m)=\frac{1}{2-e^m}$ can  be interpreted in a graph theoretic context using the idea of yokes and chains. What could be done is to take up the argument to the general generating function  $P^k(m)=\frac{e^{km}}{2-e^{m}}$ of Nelsen and Schmidt. We suspect when you glue $k$ of these yoke-chain graphs the generating function for the number of yoke-chains is related to the Nelsen-Schmidt generating function  $P^k(m)=\frac{e^{km}}{2-e^{m}}$.

\section{Acknowledgements}
Both authors acknowledge the support from Rhodes University. The first author would like also to acknowledge financial support from the DAAD-NRF scholarship of South Africa, the Levenstein Bursary of Rhodes University and the NRF-Innovation doctoral scholarship of South Africa. 

\end{document}